\documentclass[11pt,a4paper]{article}
\usepackage[margin=1in]{geometry}
\usepackage{latexsym, amsfonts, amsmath,tikz}
\usepgflibrary{arrows}
\newcommand{\be}{\begin{equation}}
\newcommand{\ee}{\end{equation}}
\newcommand{\bea}{\begin{eqnarray}}
\newcommand{\eea}{\end{eqnarray}}
\newcommand{\bean}{\begin{eqnarray*}} 
\newcommand{\eean}{\end{eqnarray*}}

\newcommand{\brray}{\begin{array}}
\newcommand{\erray}{\end{array}}
\newcommand{\ben}{\begin{equation}{nonumber}}
\newcommand{\een}{\end{equation}{nonumber}}

\newtheorem{dfn}{Definition}[section]
\newtheorem{thm}[dfn]{Theorem}
\newtheorem{lmma}[dfn]{Lemma}
\newtheorem{ppsn}[dfn]{Proposition}
\newtheorem{crlre}[dfn]{Corollary}
\newtheorem{xmpl}[dfn]{Example}
\newtheorem{rmrk}[dfn]{Remark}
\newcommand{\bdfn}{\begin{dfn}}
\newcommand{\bthm}{\begin{thm}}
\newcommand{\blmma}{\begin{lmma}}
\newcommand{\bppsn}{\begin{ppsn}}
\newcommand{\bcrlre}{\begin{crlre}}
\newcommand{\bxmpl}{\begin{xmpl}}
\newcommand{\brmrk}{\begin{rmrk}}
\newcommand{\edfn}{\end{dfn}}
\newcommand{\ethm}{\end{thm}}
\newcommand{\elmma}{\end{lmma}}
\newcommand{\eppsn}{\end{ppsn}}
\newcommand{\ecrlre}{\end{crlre}}
\newcommand{\exmpl}{\end{xmpl}}
\newcommand{\ermrk}{\end{rmrk}}

\newcommand{\cla}{{\cal A}}
\newcommand{\clb}{{\cal B}}

\newcommand{\clg}{{\cal G}}

\newcommand{\clj}{{\cal J}}
\newcommand{\clk}{{\cal K}}

\newcommand{\cln}{{\cal N}}
\newcommand{\clo}{{\cal O}}
\newcommand{\clp}{{\cal P}}
\newcommand{\clq}{{\cal Q}}

\newcommand{\clt}{{\cal T}}

\def\a*{{\cal A}_{h,*}}
\def\B{{\cal B}(h)}
\def\B1{{\cal B}_1(h)}
\def\b{{\cal B}^{\rm s.a.}(h)}
\def\b1{{\cal B}^{\rm s.a.}_1(h)}

\newcommand{\ot}{\otimes}

\newcommand{\raro}{\rightarrow}

\def \qed {$\Box$}

\def\a*{{\cal A}_{h,*}}
\def\B{{\cal B}(h)}
\def\B1{{\cal B}_1(h)}
\def\b{{\cal B}^{\rm s.a.}(h)}
\def\b1{{\cal B}^{\rm s.a.}_1(h)}

\begin{document}
\begin{center}
{\Large{\bf An example of explicit dependence of quantum symmetry on  KMS states }}\\ 
{\large {\bf Soumalya Joardar \footnote{Acknowledges support from DST INSPIRE faculty grant}   and \bf Arnab Mandal}}\\ 
\end{center}
\begin{abstract}
	We compute all the quantum symmetries of a graph with $n$-disjoint loops at the critical inverse temperature in the sense of \cite{sou_arn2} . We show that the set of non isomorphic CQG's appearing as quantum symmetry at the critical inverse temperature has a one to one correspondence with the cardinality of the set of partitions of the number of vertices.   
	\end{abstract}
{\bf Subject classification :} 46L30, 46L89\\  
{\bf Keywords:} Compact quantum group, quantum symmetry, graph $C^{\ast}$-algebra, KMS-state.
\section{Introduction}
With the advent of noncommutative geometry {\it ala} Connes, the role of compact quantum groups as symmetry (quantum) object in the realm of noncommutative geometry is now well and truly established. Since the pioneering paper of S. Wang (\cite{Wang}), where he defined the notion of quantum symmetry for finite spaces, the study of compact quantum groups as quantum symmetry object has come a long way. In \cite{Wang}, it was shown that for an $n$-point space, the quantum symmetry is remarkably large (in fact infinite dimensional as a vector space for $n\geq 4$). This drew a lot of attention of Mathematicians and the quantum symmetry concept was extended to various finite, continuous and geometric structures. For finite graphs, the concept of quantum symmetry has been studied by T. Banica, J. Bichon and others. Later the notion of quantum symmetry of Riemannian manifolds and in general spectral triple was formulated by D. Goswami and his collaborators (see \cite{Debashish},\cite{Laplace}).\\
\indent If one looks at $C^{\ast}$-algebras as noncommutative topological spaces, it is interesting to study the quantum symmetry of various concrete $C^{\ast}$-algebras. One such very well studied $C^{\ast}$-algebra is graph $C^{\ast}$-algebra associated to finite graph. It is a $C^{\ast}$-algebra generated by projections and partial isometries satisfying some relations depending on a graph. For a finite graph, the graph $C^{\ast}$-algebra is unital and one can consider it to be a noncommutative compact topological space. Recently quantum symmetry of graph $C^{\ast}$-algebras has been studied by S. Schmidt, M. Weber (see \cite{Web}) and authors of this article (see \cite{sou_arn1}, \cite{sou_arn2}). The main idea which has been used to define quantum symmetry of various structures is to identify the classical symmetry of the structure as universal object in some category of compact groups. Then one enlarges the category by replacing the compact groups by compact quantum groups (CQG for short), group morphisms by CQG morphisms and define the quantum symmetry object as the universal object of the enlarged category, if it exists. Now in general, universal object might fail to exist in the enlarged category as already observed in \cite{Wang} for matrix algebras. S. Wang remedied this problem by shrinking the category a little and demanding that the CQG action on matrix algebra to preserve the matrix trace as well. Since then it has been natural to add some kind of `volume' preserving condition on the given action.  Following this philosophy the authors of the present paper studied quantum symmetry of graph $C^{\ast}$-algebra preserving KMS states of the graph $C^{\ast}$-algebra.  In general the choice of KMS states at the critical inverse temperature is far from unique and thus the quantum symmetry depends heavily on the choice of KMS states.\\
\indent The aim of this paper is to give an example of a graph $C^{\ast}$-algebra where the set of all KMS states at critical inverse temperature form a probability simplex and compute all the quantum symmetry depending on the choice of KMS states of the graph $C^{\ast}$-algebra. In this example a particular coloring of the vertex (or equivalently a choice of partition of the number of the vertices) gives rise to a class of KMS states. Interestingly, the quantum symmetry only depends on the particular coloring of the vertices or, equivalently a particular partition of the number of the vertex. Thus, we show that only a finitely many non isomorphic CQG's arise as quantum symmetry of the particular graph $C^{\ast}$-algebra. 
\section{Preliminaries}
\subsection{Graph $C^{\ast}$ algebras and KMS states}
\label{KMSgraph}
A {\bf finite} directed graph is a collection of finitely many edges and vertices. If we denote the edge set of a graph $\Gamma$ by $E=(e_{1},\ldots,e_{n})$ and the vertex set of $\Gamma$ by $V=(v_{1},\ldots,v_{m})$ then recall the maps $s,t:E\raro V$ and the vertex matrix $D$ which is an $m\times m$ matrix whose $ij$-th entry is $k$ if there are $k$-number of edges from $v_{i}$ to $v_{j}$. 
\bdfn
$\Gamma$ is said to be without sink if the map $s:E\raro V$ is surjective. 
\edfn
\brmrk
Note that the graph $C^{\ast}$-algebra corresponding to a graph without sink is a Cuntz-Krieger algebra. Reader might see \cite{Cuntz-Krieger} for more details on Cuntz-Krieger algebra.
\ermrk
Now we recall some basic facts about graph $C^{\ast}$-algebras. Reader might consult \cite{Pask} for details on graph $C^{\ast}$-algebras. Let $\Gamma=\{E=(e_{1},...,e_{n}),V=(v_{1},...,v_{m})\}$ be a finite, directed graph without sink. We assign partial isometries $S_{i}$'s to edges $e_{i}$ for all $i=1,...,n$ and projections $p_{v_{i}}$ to the vertices $v_{i}$ for all $i=1,...,m$.
\bdfn
\label{Graph}
The graph $C^{\ast}$-algebra $C^{\ast}(\Gamma)$ is defined as the universal $C^{\ast}$-algebra generated by partial isometries $\{S_{i}\}_{i=1,\ldots,n}$ and mutually orthogonal projections $\{p_{v_{k}}\}_{k=1,\ldots,m}$ satisfying the following relations:
\begin{displaymath}
S_{i}^{\ast}S_{i}=p_{t(e_{i})}, \sum_{s(e_{j})=v_{l}}S_{j}S_{j}^{\ast}=p_{v_{l}}.
\end{displaymath}
\edfn
In a graph $C^{\ast}$-algebra $C^{\ast}(\Gamma)$, we have the following (see Subsection 2.1 of \cite{Pask}):\\
1. $\sum_{k=1}^{m}p_{v_{k}}=1$ and $S_{i}^{\ast}S_{j}=0$ for $i\neq j$.
\vspace{0.05in}\\
2. $S_{\mu}=S_{1}S_{2}\ldots S_{l}$ is non zero if and only if $\mu=e_{1}e_{2}\ldots e_{l}$ is a path i.e. $t(e_{i})=s(e_{i+1})$ for $i=1,\ldots,(l-1)$.
\vspace{0.05in}\\
3. $C^{\ast}(\Gamma)={\overline{\rm Sp}}\{S_{\mu}S_{\nu}^{\ast}:t(\mu)=t(\nu)\}$. \vspace{0.1in}\\
\indent Now we shall briefly discuss KMS-states on graph $C^{\ast}$-algebras coming from graphs without sink. For that we recall Toeplitz algebra $\clt C^{\ast}(\Gamma)$. Readers are referred to \cite{Laca} for details. Our convention though is opposite to that of \cite{Laca} in the sense that we interchange source projections and target projections. Also we shall modify the results of \cite{Laca} according to our need. Suppose that $\Gamma$ is a directed graph as before. A Toeplitz-Cuntz-Krieger $\Gamma$ family consists of mutually orthogonal projections $\{p_{v_{i}}:v_{i}\in V\}$ and partial isometries $\{S_{i}:e_{i}\in E\}$ such that $\{S_{i}^{\ast}S_{i}=p_{t(e_{i})}\}$ and
\begin{displaymath}
p_{v_{l}}\geq \sum_{s(e_{i})=v_{l}}S_{i}S_{i}^{\ast}.
\end{displaymath}
Toeplitz algebra $\clt C^{\ast}(\Gamma)$ is defined to be the universal $C^{\ast}$-algebra generated by the Toeplitz-Cuntz-Krieger $\Gamma$ family. It is clear from the definition that $C^{\ast}(\Gamma)$ is the quotient of $\clt C^{\ast}(\Gamma)$ by the ideal $\clj$ generated by
\begin{displaymath}
P:=\{p_{v_{l}}-\sum_{s(e_{i})=v_{l}}S_{i}S_{i}^{\ast}\}.
\end{displaymath}
The standard arguments give $\clt C^{\ast}(\Gamma)=\overline{\rm Sp} \{S_{\mu}S_{\nu}^{\ast}:t(\mu)=t(\nu)\}$. $\clt C^{\ast}(\Gamma)$ admits the usual gauge action $\gamma$ of $\mathbb{T}$ which descends to the usual gauge action on $C^{\ast}(\Gamma)$ given on the generators by $\gamma_{z}(S_{\mu}S_{\nu}^{\ast})=z^{(|\mu|-|\nu|)}S_{\mu}S_{\nu}^{\ast}$. Consequently it has a dynamics $\alpha:\mathbb{R}\raro {\rm Aut} \ C^{\ast}(\Gamma)$ which is lifted from $\gamma$ via the map $t\raro e^{it}$. We recall the following from \cite{Laca} (Proposition 2.1).
\bppsn
\label{exist_KMS}
Let $\Gamma$ be a finite, directed, connected graph without sink and $\gamma:\mathbb{T}\raro {\rm Aut} \ \clt C^{\ast}(\Gamma)$ be the gauge action with the corresponding dynamics $\alpha:\mathbb{R}\raro {\rm Aut} \ \clt C^{\ast}(\Gamma)$. Let $\beta\in\mathbb{R}$.\\
(a) A state $\tau$ is a ${\rm KMS}_{\beta}$ state of $(\clt C^{\ast}(\Gamma),\alpha)$ if and only if
\begin{displaymath}
\tau(S_{\mu}S_{\nu}^{\ast})=\delta_{\mu,\nu}e^{-\beta|\mu|}\tau(p_{t(\mu)}).
\end{displaymath}
(b) Suppose that $\tau$ is a ${\rm KMS}_{\beta}$ state of $(\clt C^{\ast}(\Gamma),\alpha)$, and define $\cln^{\tau}=(\cln^{\tau}_{i})\in[0,\infty)^{m}$ by $\cln^{\tau}_{i}=\phi(p_{v_{i}})$. Then $\cln^{\tau}$ is a probability measure on $V$ satisfying the subinvariance condition $D\cln^{\tau}\leq e^{\beta}\cln^{\tau}$.\\
(c) A ${\rm KMS}_{\beta}$ state factors through $C^{\ast}(\Gamma)$ if and only if $(D\cln^{\tau})_{i}=e^{\beta}\cln^{\tau}_{i}$ for all $i=1,\ldots, m$ i.e. $\cln^{\tau}$ is an eigen vector of $D$ with eigen value $e^{\beta}$.
\eppsn
\subsubsection{KMS states at critical inverse temperature}
\label{KMS}
In this subsection we state a result on existence of KMS states at inverse critical temperature on graph $C^{\ast}$-algebras coming from graphs {\bf without} sink. For that we continue to assume $\Gamma$ to be a finite, connected graph without sink and with vertex matrix $D$. We denote the spectral radius of $D$ by $\rho(D)$. With these notations, Combining Proposition 4.1 and Corollary 4.2 of \cite{Laca}, we have the following
\bppsn
\label{exist1}
The graph $C^{\ast}$-algebra $C^{\ast}(\Gamma)$ has a ${\rm KMS}_{{\rm ln}(\rho(D))}$ state say $\tau$ if and only if $\rho(D)$ is an eigen value of $D$ such that it has eigen vector with all entries non negative. If $p_{1},\ldots,p_{n}$ are the projections on the vertices, then $(\tau(p_{1}),\ldots, \tau(p_{n}))$ is an eigen vector corresponding to the eigen value $\rho(D)$ of the vertex matrix. 
\eppsn 
\subsection{Compact quantum groups as quantum automorphism groups}
\label{qaut}
In this subsection we recall the basics of compact quantum groups and their actions on $C^{\ast}$-algebras. The facts collected in this Subsection are well known and we refer the readers to \cite{Van}, \cite{Woro}, \cite{Wang} for details. 
\bdfn
A compact quantum group (CQG) is a pair $(\clq,\Delta)$ such that $\clq$ is a unital $C^{\ast}$-algebra and $\Delta:\clq\raro\clq\ot\clq$ is a
$C^{\ast}$-homomorphism satisfying\\
(i) $({\rm id}\ot\Delta)\circ\Delta=(\Delta\ot{\rm id})\circ\Delta$.\\
(ii) {\rm Sp}$\{\Delta(\clq)(1\ot\clq)\}$ and {\rm Sp}$\{\Delta(\clq)(\clq\ot 1)\}$ are dense in $\clq\ot\clq$.
\edfn
Given a CQG $\clq$, there is a canonical dense Hopf $\ast$-algebra $\clq_{0}$ in $\clq$. We shall denote the antipode by $\kappa$ and counit by $\epsilon$. Given two CQG's $(\clq_{1},\Delta_{1})$ and $(\clq_{2},\Delta_{2})$, a CQG morphism between them is a $C^{\ast}$-homomorphism $\pi:\clq_{1}\raro\clq_{2}$ such that $(\pi\ot\pi)\circ\Delta_{1}=\Delta_{2}\circ\pi$.
\bxmpl
\rm{ $S_{n}^{+}$ be the universal $C^{\ast}$-algebra generated by $\{u_{ij}\}_{i,j=1,\ldots,n}$ satisfying the following relations (see Theorem 3.1 of \cite{Wang}):
	\begin{displaymath}
	u_{ij}^{2}=u_{ij},u_{ij}^{\ast}=u_{ij},\sum_{k=1}^{n}u_{ik}=\sum_{k=1}^{n}u_{ki}=1, i,j=1,\ldots,n.
	\end{displaymath}
	It has a CQG structure with the coproduct on the generators is given by $\Delta(u_{ij})=\sum_{k=1}^{n}u_{ik}\ot u_{kj}$.}\exmpl
\bxmpl \rm{We discuss a particular unitary easy CQG. For details on unitary easy compact quantum groups the reader might consult \cite{easy}. The unitary easy CQG $H_{n}^{\infty +}$ is defined to be the universal $C^{\ast}$-algebra generated by $\{q_{ij}:i,j=1,...,n\}$ such that\vspace{0.1in}\\
	(a) the matrices $((q_{ij}))$ and $((q_{ij}^{\ast}))$ are unitaries,\vspace{0.1in}\\
	(b) all $q_{ij}$'s are normal partial isometries.\vspace{0.1in}\\
	The coproduct $\Delta$ is given on the generators by $\Delta(q_{ij})=\sum_{k=1}^{n}q_{ik}\ot q_{kj}$.
}\exmpl
\subsection{Quantum symmetry of graph $C^{\ast}$-algebras at critical inverse temperature}
\bdfn
Given a (unital) $C^{\ast}$-algebra $\cla$, a CQG $(\clq,\Delta)$ is said to act faithfully on $\cla$ if there is a unital $C^{\ast}$-homomorphism $\alpha:\cla\raro\cla\ot\clq$ satisfying\\
(i) $(\alpha\ot {\rm id})\circ\alpha=({\rm id}\ot \Delta)\circ\alpha$.\\
(ii) {\rm Sp}$\{\alpha(\cla)(1\ot\clq)\}$ is dense in $\cla\ot\clq$.\\
(iii) The $\ast$-algebra generated by the set  $\{(\omega\ot{\rm id})\alpha(\cla): \omega\in\cla^{\ast}\}$ is norm-dense in $\clq$.
\edfn
For a faithful action of a CQG $(\clq,\Delta)$ on a unital $C^{\ast}$-algebra $\cla$, there is a norm dense $\ast$-subalgebra $\cla_{0}$ of $\cla$ such that the canonical Hopf-algebra $\clq_{0}$ coacts on $\cla_{0}$.
\bdfn
(Def 2.1 of \cite{Bichon})
Given a unital $C^{\ast}$-algebra $\cla$, quantum automorphism group of $\cla$ is a CQG $(\clq,\Delta)$ acting faithfully on $\cla$ satisfying the following universal property:\\
\indent If $(\clb,\Delta_{\clb})$ is any CQG acting faithfully on $\cla$, there is a surjective CQG morphism $\pi:\clq\raro\clb$ such that $({\rm id}\ot \pi)\circ\alpha=\beta$, where $\beta:\cla\raro\cla\ot\clb$ is the corresponding action of $(\clb,\Delta_{\clb})$ on $\cla$ and $\alpha$ is the action of $(\clq,\Delta)$ on $\cla$.
\edfn
\brmrk
\label{rem1}
In general universal object might fail to exist in the above category. For the existence of universal object one generally restricts the category. In addition to $\alpha$ being a faithful action it is common to add some kind of `volume' preserving condition i.e. one assumes some linear functional $\tau$ on $\cla$ such that $(\tau\ot {\rm id})\circ\alpha(a)=\tau(a).1$ for all $a$ (see \cite{Wang}, \cite{Debashish}).
\ermrk
\bxmpl
{\rm If we take a space of $n$ points $X_{n}$ then the quantum automorphism group of the $C^{\ast}$-algebra $C(X_{n})$ is given by the CQG $S^{+}_{n}$ (see \cite{Wang}).}
\exmpl
As mentioned in Remark \ref{rem1},  for a finite, connected graph $\Gamma=(V,E)$ without sink, the category of CQG $(\clq,\Delta)$ coacting on the $C^{\ast}$-algebra $C^{\ast}(\Gamma)$ does not admit a universal object. But if we assume that the coaction to preserve a KMS state $\tau$ at inverse critical inverse temperature $\beta$ such that the eigen vector of $e^{\beta}$ has all its entries strictly positive, then universal object exists in the smaller category (see \cite{sou_arn2}).We call such a universal object which depends on the state $\tau$, the quantum automorphism group of $C^{\ast}(\Gamma)$ at the state $\tau$ and denote it by $Q_{aut}(C^{\ast}(\Gamma))_{\tau}$.
\brmrk
In case of strongly connected graphs, the corresponding graph $C^{\ast}$-algebra has a unique KMS state at its inverse critical temperature and we call the quantum automorphism group preserving that state to be {\bf the} quantum automorphism group of the graph $C^{\ast}$-algebra at the critical inverse temperature (see \cite{sou_arn2}). But in general there can be infinitely many KMS states at the inverse critical temperature and thus apriori many quantum automorphism groups might exist depending on the choice of states.
\ermrk
In general quantum symmetry of graph $C^{\ast}$-algebras are not of Kac type. Nevertheless, the antipode has the following special property (see the proof of Lemma 3.4 of \cite{sou_arn2}) which we shall use later in the paper . In the following Lemma, $\Gamma$ is a finite, connected graph with $n$-edges so that $C^{\ast}(\Gamma)$ is generated by $n$ partial isometries $S_{i}$ for $i=1,\ldots,n$ and $\tau$ be a KMS state at the critical inverse temperature.
\blmma
\label{kappa}
Let $\alpha$ be the coaction of $Q_{aut}(C^{\ast}(\Gamma))_{\tau}$ given by $\alpha(S_{i})=\sum_{j=1}^{n} S_{j}\ot q_{ji}, i=1,\ldots,n$. Then $\kappa(q_{ij})=\lambda_{ij}q_{ji}^{\ast}$ for some constants $\lambda_{ij}$'s for all $i,j=1,\ldots,n$.
\elmma
\subsection{ Free  wreath product of CQG}
	In this paper, we shall use the following universal property of the free product of CQG's (see \cite{Wangfree})
	\bppsn
	Given an arbitrary family of CQG's $(\clq_{\lambda},\Delta_{\lambda})$ and a CQG $(\clq,\Delta)$ along with CQG-morphisms $\alpha_{\lambda}:\clq_{\lambda}\raro\clq$, there is a CQG morphism $\alpha:\star_{\lambda}\clq_{\lambda}\raro\clq$ such that $\alpha_{\lambda}=\alpha\circ i_{\lambda}$ where $i_{\lambda}:\clq_{\lambda}\raro\star_{\lambda}\clq_{\lambda}$ are natural inclusions for all $\lambda$, where $\star$ stands for the free product of CQG's.
	\eppsn
\bdfn (see \cite{Bichon2})
Let $(\clq,\Delta)$ be a CQG and $n\in\mathbb{N}$. Then the free wreath product of $(\clq,\Delta)$ with the quantum permutation group $S_{n}^{+}$, as a $C^{\ast}$-algebra, is the quotient of the $C^{\ast}$-algebra $(\underbrace{\clq\star...\star\clq}_{n-copies}\star S_{n}^{+})$ by the two sided ideal generated by the elements
\begin{displaymath}
	i_{k}(q)x_{ki}-x_{ki}i_{k}(q),1\leq i,k\leq n, q\in\clq,
	\end{displaymath}
where $i_{k}$ is the inclusion of k-th copy of $\clq$ in $\underbrace{\clq\star...\star\clq}_{n-copies}$  and $\{x_{ki};1\leq i,k\leq n\}$ are generators of $S_{n}^{+}$. The free wreath product CQG is denoted by $\clq\wr S_{n}^{+}$.  The coproduct on the generators is given by
\begin{displaymath}
\Delta(x_{ij})=\sum x_{ik}\ot x_{kj}; \Delta(i_{k}(q))=\sum_{j=1}^{n}i_{k}\ot i_{j}(\Delta(q))(x_{kj}\ot 1).
\end{displaymath} 
\edfn
We recall the following (see Theorem 3.4 of \cite{fusion})
\bppsn
\label{wreath}
$H^{\infty+}_{n}\cong C(S^{1})\wr S_{n}^{+}$.
\eppsn
\section{Main section}
In this section we concentrate on quantum symmetry of a particular graph $\clg=(V,E)$ consisting of $n$-disjoint loops. It is easy to see that the vertex matrix is given by ${\rm Id}_{n\times n}$ (see \cite{sou_arn2}).  The graph $C^{\ast}$-algebra $C^{\ast}(\clg)$ is generated by $n$-partial isometries $S_{i}$ and mutually orthogonal projections $\{p_{i}:i=1,\ldots,n\}$ such that 
\begin{displaymath}
S_{i}^{\ast}S_{i}=S_{i}S_{i}^{\ast}=p_{i}.
\end{displaymath}
\blmma
\label{word}
For $C^{\ast}(\clg)$, for any $n\in\mathbb{N}$, $S_{i_{1}}\ldots S_{i_{n}}S_{j_{n}}^{\ast}\ldots S_{j_{1}}^{\ast}=\delta_{i_{1},\ldots,i_{n},j_{n},\ldots, j_{1}}p_{i_{1}}$.
\elmma 
{\it Proof}:\\
It follows from the simple observation that for $i\neq j$, $S_{i}S_{j}=S_{i}S_{j}^{\ast}=S_{i}^{\ast}S_{j}=0$ and for $i=j$, $S_{i}^{\ast}S_{i}=S_{i}S_{i}^{\ast}=p_{i}$.\qed\vspace{0.1in}\\
Since the eigen space of ${\rm Id}_{n\times n}$ is the whole of $\mathbb{R}^{n}$, by Proposition \ref{exist1}, the KMS states of $C^{\ast}(\clg)$ at the critical temperature $0$ are given by 
\begin{displaymath}
\tau(S_{i}S_{i}^{\ast})=c_{i},\sum_{i}c_{i}=1, 0<c_{i}<1,
\end{displaymath} where $(c_{1},\ldots,c_{n})$ is an eigen vector of ${\rm Id}_{n\times n}$. Given a partition $\clp=(m_{1},\ldots, m_{k})$ of $n$,   we can number the edges so that $C^{\ast}(\clg)$ is generated by the partial isometries $\{S^{(i)}_{\mu}:\mu=1,\ldots, m_{i} \ and \ i=1,\ldots,k\}$. We denote the projections $S^{(i)\ast}_{\mu}S^{(i)}_{\mu}(=S^{(i)}_{\mu}S^{(i)\ast}_{\mu})$ by $p^{(i)}_{\mu}$. Corresponding to the partition, we have a class of KMS states denoted by $[\clp]$.  $\tau\in[\clp]$ is given on the projections by 
\begin{displaymath}
\tau(p^{(i)}_{\mu})=c^{(i)}_{\mu},\mu=1,\ldots,m_{i}:i=1,\ldots,k,
\end{displaymath}
so that $c^{(i)}_{\mu}=c^{(i)}_{\nu}$ for all $\mu,\nu=1,\ldots,m_{i}$,  $c^{(i)}_{\mu}\neq c^{(j)}_{\nu}$ for any $i\neq j$ and $\sum_{i=1,\mu=1}^{i=k,\mu=m_{i}}c^{(i)}_{\mu}=1$. Varying the partition $\clp$, we can get all the KMS states of the graph $C^{\ast}$-algebra $C^{\ast}(\clg)$. More precisely if we denote the set of KMS states of $C^{\ast}(\clg)$ at the temperature $0$ by $\clk_{\clg}$, then $\clk_{\clg}=\{\tau\in[\clp]:\clp \ is \ a \ partition \ of \ n\}$. Now we can state the main Theorem of  the paper.
\bthm
\label{main theorem}
Let $\clg$ be the graph consisting of $n$-disjoint loops and $\clp=(m_{1},\ldots,m_{k})$ be a partition of $n$. Then for any KMS state $\tau\in[\clp]$ of $C^{\ast}(\clg)$ at the critical inverse temperature $0$, $Q_{aut}(C^{\ast}(\clg))_{\tau}\cong  (\star_{i=1}^{k}(C(S^{1})\wr S_{m_{i}}^{+}))$.
\ethm
Given a partition $\clp=(m_{1},\ldots,m_{k})$,  $Q_{aut}(C^{\ast}(\clg))_{\tau}$ is generated algebraically by the elements $q^{(ji)}_{\nu\mu}$ where the corresponding action on $C^{\ast}(\clg)$ is given by 
\begin{eqnarray}
\label{action}
\alpha(S^{(i)}_{\mu})=\sum_{j,\nu}S^{(j)}_{\nu}\ot q^{(ji)}_{\nu\mu}.
\end{eqnarray} 
\blmma
\label{partial isometries}
$q^{(ji)}_{\nu\mu}$'s are partial isometries for all $i,j=1,\ldots,k$ and $\mu=1,\ldots,m_{i}$ and $\nu=1,\ldots,m_{j}$. 
\elmma
{\it Proof}:\\
Using the fact that  $S^{(i)\ast}_{\mu}S^{(j)}_{\nu}=S^{(i)}_{\mu}S^{(j)\ast}_{\nu}=\delta_{ij}\delta_{\mu\nu}p^{(i)}_{\mu}$, we get
\begin{displaymath}
\alpha(S^{(i)}_{\mu})=\alpha(S^{(i)}_{\mu}S^{(i)\ast}_{\mu}S^{(i)}_{\mu})=\sum S^{(j)}_{\nu}S^{(j)\ast}_{\nu}S^{(i)}_{\nu}\ot q^{(ji)}_{\nu\mu}q^{(ji)\ast}_{\nu\mu}q^{(ji)}_{\nu\mu}=\sum S^{(j)}_{\nu}\ot q^{(ji)}_{\nu\mu}q^{(ji)\ast}_{\nu\mu}q^{(ji)}_{\nu\mu}.
\end{displaymath}
But $\alpha(S_{\mu}^{(i)})=\sum S^{(j)}_{\nu}\ot q^{(ji)}_{\nu\mu}$. Hence comparing coefficients, we get $q^{(ji)}_{\nu\mu}q^{(ji)\ast}_{\nu\mu}q^{(ji)}_{\nu\mu}=q^{(ji)}_{\nu\mu}$ for all $i,j,\mu,\nu$.\qed
\blmma
\label{q_{ij}*q_{ik}}
$q^{(il)\ast}_{\mu\xi}q^{(jl)}_{\nu\xi}=0$ whenever $i\neq j$ or $\mu\neq\nu$.
\elmma
{\it Proof}:\\
Using the fact that if $i\neq j$ or $\mu\neq\nu$, $S^{(i)\ast}_{\mu}S^{(j)}_{\nu}=0$, we get for $i\neq j$ or $\mu\neq\nu$,
\begin{displaymath}
0=\alpha(S^{(i)\ast}_{\mu}S^{(j)}_{\nu})=\sum S^{(l)\ast}_{\xi}S^{(l)}_{\xi}\ot q^{(li)\ast}_{\xi\mu}q^{(lj)}_{(\xi\nu)}.
\end{displaymath}
Using the linear independence of the projections $S^{(l)\ast}_{\xi}S^{(l)}_{\xi}$, we get $q^{(li)\ast}_{\xi\mu}q^{(lj)}_{\xi\nu}=0$. Applying $\kappa$ to this identity and using Lemma \ref{kappa}, we finish the proof of the Lemma.
 \qed
\blmma
\label{block}
$q^{(ji)}_{\xi\mu}=0$ whenever $i\neq j$.
\elmma
{\it Proof}:\\
Let $i\neq j$. Fix some $\mu\in\{1,\ldots,m_{i}\}$. Then $\alpha(S^{(i)}_{\mu}S_{\mu}^{(i)\ast})=\sum S^{(l)}_{\nu}S^{(l)\ast}_{\nu}\ot q^{(li)}_{\nu\mu}q^{(li)\ast}_{\nu\mu}$. Applying $(\tau\ot{\rm id})$ to the both side of this equation and using the $\tau$-preserving property of $\alpha$, we get $\tau(p^{(i)}_{\mu})=\sum\tau(p^{(l)}_{\nu})q^{(li)}_{\nu\mu}q^{(li)\ast}_{\nu\mu}$. Multiplying both sides by $q^{(ji)}_{\xi\mu}$ and using Lemma \ref{q_{ij}*q_{ik}} and Lemma \ref{partial isometries}, we get
\begin{displaymath}
\tau(p^{(i)}_{\mu})q^{(ji)}_{\xi\mu}=\tau(p^{(j)}_{\xi})q^{(ji)}_{\xi\mu}.
\end{displaymath}
But by definition of $\tau$, for $i\neq j$, $\tau(p^{(i)}_{\mu})\neq \tau(p^{(j)}_{\xi})$. Hence $q^{(ji)}_{\xi\mu}=0$ whenever $i\neq j$.\qed.\vspace{0.1in}\\
{\it Proof of Theorem \ref{main theorem}}:\\
Let $\alpha$ be the action of $Q_{aut}(C^{\ast}(\clg))_{\tau}$ on $C^{\ast}(\clg)$ as given in (\ref{action}). Then by virtue of the Lemma \ref{block}, we can rewrite the action as 
\begin{displaymath}
\alpha(S_{\mu}^{(i)})=\sum_{\nu=1}^{m_{i}}S_{\nu}^{(i)}\ot q_{\nu\mu}^{(i)},\mu=1,\ldots,m_{i} \ and \ i=1,\ldots,k,
\end{displaymath}
where $\tau(S_{\mu}^{(i)\ast}S_{\mu}^{(i)})=\tau(S_{\nu}^{(i)\ast}S_{\nu}^{(i)})$ for $\mu,\nu=1,\ldots, m_{i}$ and $i=1,\ldots,k$ with $\tau(S_{\mu}^{(i)\ast}S_{\mu}^{(i)})\neq\tau(S_{\mu}^{(j)\ast}S_{\mu}^{(j)})$ for $i\neq j$. Since $\alpha$ preserves $\tau$, the matrices $((q_{\mu\nu}^{(i)})), ((q_{\mu\nu}^{(i)\ast}))$ are $m_{i}\times m_{i}$ unitary matrices for all $i=1,\ldots,k$. Also using the same arguments as in Example 2 of \cite{sou_arn2}, we can show that each $q_{\mu\nu}^{(i)}$ is normal partial isometries. So by universal property, there is a surjective CQG morphism from $H^{\infty +}_{m_{i}}\overset{\mathrm Prop. \ \ref{wreath}}\cong (C(S^{1})\wr S_{m_{i}}^{+})$ onto $Q_{aut}(C^{\ast}(\clg))_{\tau}$ for all $i=1,\ldots,k$. So by the universal property of free product, there is  a surjective CQG morphism say $\Lambda$ from $\star_{i=1}^{k}(C(S^{1})\wr S_{m_{i}}^{+})$ onto $Q_{aut}(C^{\ast}(\clg))_{\tau}$.\\
\indent For the inverse CQG morphism, let us denote the generators of $(C(S^{1})\wr S_{m_{i}}^{+})$ by $\{q^{(i)}_{\mu\nu}:\mu,\nu=1,\ldots,m_{i}\}$ for $i=1,\ldots,k$.  We define $\alpha:C^{\ast}(\clg)\raro C^{\ast}(\clg)\ot (\star_{i=1}^{k}(C(S^{1})\wr S_{m_{i}}^{+}))$ by
\begin{eqnarray}
\label{alpha}
	\alpha(S^{(i)}_{\mu})=\sum_{\nu=1}^{m_{i}}S_{\nu}^{(i)}\ot q_{\nu\mu}^{(i)},i=1,\ldots,k \ and \ \mu=1,\ldots,m_{i}.
	\end{eqnarray}
Using the fact that $q_{\nu\mu}^{(i)}$'s are normal elements and $S_{\nu}^{(i)\ast}S_{\xi}^{(i)}=S_{\nu}^{(i)}S_{\xi}^{(i)\ast}=0$ for $\xi\neq\nu$, we get that 
\begin{displaymath}
\alpha(S_{\mu}^{(i)})^{\ast}\alpha(S_{\mu}^{(i)})=\alpha(S_{\mu}^{(i)})\alpha(S_{\mu}^{(i)})^{\ast}=\sum_{\nu}p^{(i)}_{\nu}\ot q^{(i)\ast}_{\nu\mu}q^{(i)}_{\nu\mu}, \mu=1,\ldots,m_{i},i=1,\ldots,k.
\end{displaymath}
If we define $\alpha(p^{(i)}_{\mu})$ by $\alpha(S_{\mu}^{(i)})^{\ast}\alpha(S_{\mu}^{(i)})$, then it is easy  to see that $\alpha(p^{(i)}_{\mu})^{\ast}=\alpha(p^{(i)}_{\mu})$.
\begin{eqnarray*}
	\alpha(p^{(i)}_{\mu})\alpha(p^{(j)}_{\nu})&=&(\sum p^{(i)}_{\xi}\ot q^{(i)\ast}_{\xi\mu}q^{(i)}_{\xi\mu})(\sum p^{(j)}_{\eta}\ot q^{(j)\ast}_{\eta\nu}q^{(j)}_{\eta\nu})\\
	&=& \delta_{ij}\sum p^{(i)}_{\xi}\ot q^{(i)\ast}_{\xi\mu}q^{(i)}_{\xi\mu}q^{(i)\ast}_{\xi\nu}q^{(i)}_{\xi\nu}.
	\end{eqnarray*}
Since for a fixed $\xi$,  $q^{(i)\ast}_{\xi\mu}q^{(i)}_{\xi\mu}$'s are mutually orthogonal projections , $\alpha(p^{(i)}_{\mu})\alpha(p^{(j)}_{\nu})=\delta_{ij}\delta_{\mu\nu}\sum p^{(i)}_{\xi}\ot q^{(i)\ast}_{\xi\mu}q^{(i)}_{\xi\mu}$ i.e. $\alpha(p^{(i)}_{\mu})\alpha(p^{(j)}_{\nu})=\delta_{ij}\delta_{\mu\nu}\alpha(p^{(i)}_{\mu})$. Hence by universal property of $C^{\ast}(\clg)$, we in deed have a well defined $C^{\ast}$-homomorphism $\alpha$ given by (\ref{alpha}). That $\alpha$ is a coaction is trivial to check. In order to show that $\alpha$ preserves $\tau$,  for $r\neq s$,  by the definition of the KMS state, $\tau(S_{\mu_{1}}^{(i_{1})}\ldots S_{\mu_{r}}^{(i_{r})}S_{\nu_{s}}^{(j_{s})\ast}\ldots S_{\nu_{1}}^{(j_{1})\ast})=0$. Since $\alpha$ is linear, 
  \begin{displaymath}
  (\tau\ot{\rm id})\circ\alpha(S_{\mu_{1}}^{(i_{1})}\ldots S_{\mu_{r}}^{(i_{r})}S_{\nu_{s}}^{(j_{s})\ast}\ldots S_{\nu_{1}}^{(j_{1})\ast})=0.
  \end{displaymath}
  Hence for $r\neq s$, $(\tau\ot{\rm id})\circ\alpha(S_{\mu_{1}}^{(i_{1})}\ldots S_{\mu_{r}}^{(i_{r})}S_{\nu_{s}}^{(j_{s})\ast}\ldots S_{\nu_{1}}^{(j_{1})\ast})=\tau(S_{\mu_{1}}^{(i_{1})}\ldots S_{\mu_{r}}^{(i_{r})}S_{\nu_{s}}^{(j_{s})\ast}\ldots S_{\nu_{1}}^{(j_{1})\ast})1$. For $r=s$,
\begin{eqnarray*}
&&(\tau\ot{\rm id})\circ\alpha(S_{\mu_{1}}^{(i_{1})}\ldots S_{\mu_{r}}^{(i_{r})}S_{\nu_{r}}^{(j_{r})\ast}\ldots S_{\nu_{1}}^{(j_{1})\ast})\\
&=& \delta_{i_{1},\ldots,i_{r},j_{r},\ldots,j_{1}}\delta_{\mu_{1},\ldots,\mu_{r},\nu_{r},\ldots,\nu_{1}}(\tau\ot{\rm id})\circ\alpha (p_{\mu_{1}}^{(i_{1})}) \ (by \ Lemma \ \ref{word})\\
&=& \delta_{i_{1},\ldots,i_{r},j_{r},\ldots,j_{1}}\delta_{\mu_{1},\ldots,\mu_{r},\nu_{r},\ldots,\nu_{1}}\sum_{\nu} \tau(p^{(i_{1})}_{\nu})q^{(i_{1})}_{\nu\mu_{1}}q^{(i_{1})\ast}_{\nu\mu_{1}}\\
&=& \delta_{i_{1},\ldots,i_{r},j_{r},\ldots,j_{1}}\delta_{\mu_{1},\ldots,\mu_{r},\nu_{r},\ldots,\nu_{1}}\tau(p^{(i_{1})}_{\mu_{1}}) \ (since \ \tau(p^{(i_{1})}_{\nu})=\tau(p^{(i_{1})}_{\mu_{1}}) \ \forall \nu \ and \ \sum_{\nu}q^{(i_{1})}_{\nu\mu_{1}}q^{(i_{1})\ast}_{\nu\mu_{1}}=1),
\end{eqnarray*}
which is equal to $\tau(S_{\mu_{1}}^{(i_{1})}\ldots S_{\mu_{r}}^{(i_{r})}S_{\nu_{r}}^{(j_{r})\ast}\ldots S_{\nu_{1}}^{(j_{1})\ast})1$. Therefore for all $\mu_{1},\ldots,\mu_{r},\nu_{1},\ldots,\nu_{r}$ and $i_{1},\ldots,i_{r},j_{r},\ldots,j_{1}$, 
\begin{displaymath}(\tau\ot{\rm id})\circ\alpha(S_{\mu_{1}}^{(i_{1})}\ldots S_{\mu_{r}}^{(i_{r})}S_{\nu_{r}}^{(j_{r})\ast}\ldots S_{\nu_{1}}^{(j_{1})\ast})=\tau(S_{\mu_{1}}^{(i_{1})}\ldots S_{\mu_{r}}^{(i_{r})}S_{\nu_{r}}^{(j_{r})\ast}\ldots S_{\nu_{1}}^{(j_{1})\ast})1.\end{displaymath} Since for any finite graph $\Gamma$, $C^{\ast}(\Gamma)={\overline{\rm Sp}}\{S_{\mu}S_{\nu}^{\ast}:t(\mu)=t(\nu)\}$, by linearity and continuity of $\alpha$ and $\tau$, we conclude that $(\tau\ot{\rm id})\circ\alpha(a)=\tau(a).1$ for all $a\in C^{\ast}(\clg)$, i.e. $\alpha$ preserves $\tau$. So, by universal property of $Q_{aut}(C^{\ast}(\clg))_{\tau}$, there is a surjective CQG morphism from $Q_{aut}(C^{\ast}(\clg))_{\tau}$ onto $(\star_{i=1}^{k}(C(S^{1})\wr S_{m_{i}}^{+}))$ . It is straightforward to see that this morphism is in deed the inverse morphism of $\Lambda$. Hence for any partition $\clp=(m_{1},\ldots,m_{k})$ of $n$ and any $\tau\in[\clp]$,
\begin{displaymath} Q_{aut}(C^{\ast}(\clg))_{\tau}\cong (\star_{i=1}^{k}(C(S^{1})\wr S_{m_{i}}^{+})).
\end{displaymath}
\qed 
\brmrk
Observe that the quantum symmetry only depends on the class $[\clp]$. Specifying the class $[\clp]$ is equivalent to coloring the vertices. For example given a partition $\clp=(m_{1},\ldots,m_{k})$, we have a $k$-coloring of the vertices of the graph so that $m_{i}$ vertices gets $i$-th color. $Q_{aut}(C^{\ast}(\clg))_{\tau}$ preserves any $\tau\in[\clp]$ is equivalent to demanding that  $Q_{aut}(C^{\ast}(\clg))_{\tau}$ preserves such coloring and vice versa.
\ermrk
  
Soumalya Joardar \\
Theoretical Science Unit,\\ 
JNCASR, Bangalore-560064, India\\ 
email: soumalya.j@gmail.com 
\vspace{0.1in}\\
Arnab Mandal\\
Presidency University, College Street\\
Kolkata-700073, India\\
email: arnabmaths@gmail.com

\end{document}